\documentclass[11pt]{amsart}

\usepackage{amsmath}
\usepackage{amsfonts}
\usepackage{amssymb} 
\usepackage{amsthm} 
\usepackage{mathptmx}
\usepackage{mathrsfs}
\usepackage{latexsym}
\usepackage{times}
\usepackage[mathscr]{euscript}
\usepackage[isolatin]{inputenc}

\usepackage[pagebackref]{hyperref}
\usepackage[backrefs, alphabetic,initials]{amsrefs}

\DeclareMathAlphabet{\mathpzc}{OT1}{pzc}{m}{it}

\theoremstyle{definition}

\theoremstyle{remark}

\newcommand\Cc{\mathcal C}

\newcommand\CC{\mathbb C}
\newcommand\NN{\mathbb N}

\renewcommand\hat\widehat
\renewcommand\tilde\widetilde

\begin{document}

\font\eightrm=cmcsc10

\font\eightr=cmr10 at 14pt
\font\eight=cmr10 at 12pt

\def\R{{I\kern-0.3emR}}
\def\P{{I\kern-0.3emR}}
\def\K{{I\kern-0.3emR}}

\baselineskip=0mm
\lineskiplimit=1pt
\lineskip=2mm

\title[ON TWO THEOREMS ABOUT SYMPLECTIC REFLECTION ALGEBRAS]{ON TWO
  THEOREMS ABOUT SYMPLECTIC REFLECTION ALGEBRAS}

\author{Georges Pinczon}

\address{Institut de Mathématiques de Bourgogne, Université de
  Bourgogne, B.P. 47870, F-21078 Dijon Cedex, France}

\email{gpinczon@u-bourgogne.fr}


\subjclass[2000]{}

\date{\today}

\begin{abstract} 
We give a new proof and an improvement of two Theorems of J. Alev,
M.A. Farinati, T. Lambre and A.L. Solotar [1] : the first one about
Hochschild cohomology spaces of some twisted bimodules of the Weyl
Algebra $W$, and the second one about Hochschild cohomology spaces of
the smash product $G * W$ ($G$ a finite subgroup of $SP (2n)$) and, as
a consequence, we then give a new proof of a Theorem of P. Etingof and
V. Ginzburg [12], which shows that the Symplectic Reflection Algebras
are deformations of $G * W$ (and, in fact, all possible ones).
\end{abstract}

\date{\today}

\maketitle

\vskip15pt
{\bf INTRODUCTION}
\vskip10pt
This paper belongs to a very fascinating context of non-commutative geometry, for which we refer to [12] and [11]: in these papers, this context is perfectly described an developed, with many deep results, applications and examples (and all references). The problem is the study of deformations of $G * W,$ where $W$ is a Weyl Algebra, and $G$ a finite subgroup of $SP (2n)$. Following Gerstenhaber, one has to find $H^2 (G * W)$ (first order deformations), then $H^3 (G * W)$ (obstructions), and then universal models of deformations of $G * W.$ It turns out that these models are the Symplectic Reflection Algebras of [12]. Let us quickly explain how this program was worked out.
\vskip5pt
In [1], J. Alev, M.A. Farinati, T. Lambre and A.L. Solotar have found:

(1) (A.F.L.S. Theorem 1) the dimension of the Hochschild cohomology
spaces $H^k (W_\sigma)$ of the twisted W-bimodule $W_\sigma$ ($\sigma$
a diagonalizable element of $SP (2n)$).

(2) (A.F.L.S. Theorem 2) the dimension of the Hochschild cohomology spaces $H^k (G * W).$

Only particular cases of the A.F.L.S. Theorems were known before: when $\sigma = Id,$ i.e the Hochschild cohomology of the Weyl Algebra ([20] [10]), or when $\sigma$ is the parity of $W$ ([17]). Part 2 gives the cohomology of $W^G,$ since $W^G$ and $G * W$ are Morita-equivalent ([15].[1]). Some more information about $H^\bullet (G * W),$ including a description of the cup-product, was obtained in [3]. For applications of the A.F.L.S. Theorem see. [1] and [3]. Later, a new proof of the AFLS Theorems was given by P. Etingof [11].

In [12] (see also [11]), among many other results, P. Etingof and V. Ginzburg have completely solved the problem of deformations of $G * W$, showing that their Symplectic Reflections Algebras are non trivial (algebraic) deformations of $G * W$, and describe (up to equivalence and change of parameter) all possible deformations (the E.G. Theorem).

The E.G. Theorem belongs to non-commutative geometry (see [12]), and also to deformation quantization theory [5], since the Symplectic Reflection Algebras are natural generalizations of algebras used to quantize Calogero-Moser systems (see [11], [12]).

The goal of the present paper is to give new proofs of the
A.F.L.S. Theorems, and also of the E.G. Theorem. We do not pretend
that our proofs are simpler than the original ones, they are
different, and we believe that new different proofs of deep results
may be of interest. Moreover, we prove a significant amelioration of
the A.F.L.S. Theorems, let us call it the C.A. Theorems, which can be
used to simplify the original proof of the E.G. Theorem, and is an
essential argument to build the new proof of the E.G. Theorem.

Let us now describe the sections of the paper, and the results.

(1) In section 1, we revisit the Koszul complex of the Weyl Algebra,
as defined in [17] [1]. Let us quickly explain why. When reading
papers [17], [1], [3], [12], [11], one has the feeling that the
A.F.L.S. Theorems are not achieved: they give the dimensions of
cohomology spaces, but no information about cocycles themselves. For
instance, in [12], P. Etingof and V. Ginzburg have to prove that their
deformation is not trivial (Lemma (2.17)), and the proof is not at all
trivial. To understand what is missing, let us give another
example. Let $S$ be a polynomial algebra, the Hochschild cohomology of
$S$ reduces to the cohomology of the Koszul complex, which is infinite
dimensional at all degrees. But one can obtain much more information
if one remarks that the Koszul complex is a subcomplex of the Bar
resolution: it easily follows that the cohomology is the space of
skewsymmetric multivectors, and many other useful consequences. Now,
we come back to the Weyl Algebra $W.$ In that case, the Koszul complex
is not a subcomplex of the Bar resolution, nevertheless the following
holds:

\vskip5pt LEMMA \vskip5pt {\it The Koszul complex of $W$ is a
  subcomplex of the normalized Bar resolution.}  \vskip6pt To our
knowledge, this result has never been stated up to now. As a
consequence, we prove:

\vskip10pt THEOREM 1: \vskip5pt {\it Given \ a \ W-bimodule \ $M,$
  \ one \ has \:}

{\it (1) The restriction map, from the Hochschild complex of $M$ to the Koszul complex of $M,$ induces an isomorphism in cohomology.

(2) Let $V = span (p_i, q_i, i = 1...n),$ and $\Lambda$ the exterior algebra of $V.$ Any Koszul $k$-cocycle is the restriction to $\Lambda^k$ of a Hochschild $k$-cocycle.}

(2) In section 2, we prove Part.1 of the C.A. Theorem. Let $\sigma$ be
a diagonalizable element of $SP (2n),$ $V_\sigma:= Range (\sigma -
Id), \ 2 k_\sigma:= dim \ V_\sigma,$ $\omega_\sigma$ the form defined
by $\omega_\sigma (v_1,...,v_{2 k_\sigma}):= {1 \over {{k_\sigma}!}}
\omega^{k_\sigma} (p_\sigma (v_1), ...,$ $p_\sigma (v_{2 k
  \sigma_\sigma})),$ where $\omega$ is the canonical two-form and
$p_\sigma$ the projection on $V_\sigma$ coming from $V = V_\sigma
\ \oplus \ Im (\sigma - Id),$ $W_\sigma$ the twisted $W$-bimodule
associated to $\sigma.$

\vskip10pt 

C.A. THEOREM 1 \vskip5pt {\it (1) (A.F.L.S. Theorem 1) $dim
  \ H^{2k_\sigma} (W_\sigma) = 1,$ and $H^k (W_\sigma) = \{ 0 \},$ if
  $k \not = 2 k_\sigma.$}

{\it (2) There exists a Hochschild cocycle $\Omega_\sigma \in Z^{2 k_\sigma} (W_\sigma)$ such that:} $\Omega_\sigma \ \vert_{\Lambda^{2 k_\sigma}} = \omega_\sigma, \ {\rm and} \ H^{2k_\sigma} (W_\sigma) = \CC . \Omega_\sigma.$

{\it (3) Let $S$ be a finite subgroup of $SP (2n).$ We assume that $S$
  commutes with $\sigma,$ and denote by $H^\bullet_S (W_\sigma)$ the
  $S$-invariant Hochschild cohomology. Then $H^\bullet_S (W_\sigma) =
  H^\bullet (W_\sigma),$ and the cocycle $\Omega_\sigma$ in (2) can be
  chosen to be $S$-invariant.}

Let us say of few words about the proof. As in [17], [1], or [11], we
introduce the Koszul complex, but then we follow arguments of [17]:
we replace the differential by an equivalent one, split the new
complex, and introduce explicit homotopies to deduce (1), then (2) and
(3) follow from Theorem 1. Let us remark that the equivalence we use
is natural in this context: it is built using operators which appear
when showing that standard ordering and Weyl ordering define
equivalent star-products.

(3) In section (3), we prove Part.2 of the C.A. Theorem. This is an easy consequence of Part.1: $\CC [G]$ is separable, so $H^\bullet_{{\CC} [G]} (G * W) = H^\bullet (G * W) \ [13],$ and using a convenient description of $\CC [G]$-relative cocycles, the result follows.

Let $\Gamma$ be the set of conjugacy classes of $G$ and $\Gamma_{2k}
:= \{ \gamma \in \Gamma / k_\sigma = k, \ \forall \ \sigma \in \Gamma
\}.$ 

\vskip15pt

C.A. THEOREM 2
\vskip5pt
{\it (1) (A.F.L.S. Theorem 2) $H^k (G * W) = \{ 0 \},$ if $k$ is odd, and $dim \ H^{2k} (G * W) = card \ \Gamma_{2k}.$

(2) Assuming that $\Gamma_{2k} \not = \phi,$ let $\lambda \in \CC^{\Gamma_{2k}}.$ There exists a cocycle $C_\lambda \in Z^{2 k}_{\CC [G]} (G * W)$ such that}
\vskip0.5pt
$\ \ \ \ \ \ \ \ \ \ C_\lambda \ (X_1 \ \wedge ... \wedge \ X_{2k}) = \displaystyle {\sum_{\gamma \in \Gamma_{2k}}} \ \lambda (\gamma) \ \displaystyle {\sum_{g \in \gamma}} \ \ \omega_g (X_1,...,X_{2k}) \ \otimes \ g, \ \forall \ X_i \in V.$
\vskip0.5pt
{\it $C_\lambda$ is a coboundary if and only if $\lambda = 0,$ and the map $\lambda \rightarrow C_\lambda$ induces an isomorphism from $\CC^{\Gamma_{2k}}$ onto $H^{2k}_{ \CC \ [G]} (G * W).$ 

(3) If a cocycle $C \in Z^{2k}_{\CC [G]} \ (G * W)$ vanishes on $\Lambda^{2k},$ then $C$ is a coboundary.}

(2) and (3) are of interest since they describe the "emerged part" of cocycles. For instance Lemma (2.17) of [12] is a consequence of (2) and (3).

(4) In section (4), we give a new proof of the E.G. Theorem. Assuming
$\Gamma_2 \not = \phi,$ let $H_{\hbar \lambda}$ be the Symplectic
Reflection Algebra, with $\lambda \in \ \CC^{\Gamma_2}, \lambda \not =
0.$ This algebra is defined by generators and relations
(SRA-relations), see (4.2). We prove: 

\vskip5pt E.G. THEOREM \vskip5pt {\it There exists a non trivial
  polynomial $\CC [G]$-relative deformation of $G * W$ where the
  SRA-relations hold. The subalgebra $(G * W)[\hbar]$ of this
  deformation is isomorphic to $H_{\hbar \lambda}.$ We have $<
  X_1,..., X_k > := \ {1 \over {k!}} \ \displaystyle {\sum_{\sigma \in
      \Sigma}} \ X_{\sigma (1)} \ \displaystyle\mathop{*}_\hbar
  \oalign{\cr} \ ... \displaystyle\mathop{*}_\hbar \oalign{\cr}$
  $X_{\sigma (k)}= X_1 ... X_k, \ \forall \ X_i \in V$ (the last
  product is computed using the abelian product on $W$).}

The last properties give the P.B.W. property: given a basis $\{ e_1,..., e_{2n} \}$ of $V,$ $\{ e^{i_1}_1 ... e^{i_{2n}}_{2n} \ \otimes \ g, \ i_j \in \NN, \ g \in G \}$ is, resp.: when $\hbar$ is formal, a $\ \CC  [\hbar]$-basis of $H_{\hbar \lambda} = (G * W) \ [\hbar],$ resp.: when $\hbar \in \CC,$ a basis of $H_{\hbar \lambda}.$

As stated, the theorem is completely equivalent to the E.G.-Theorem, only the order of the claims and the proof differ. Let us give some details:

The original proof in [12] has two steps: 

\vskip5pt $\bullet$ First step (main one): $H_{\hbar \lambda},$ is a
deformation due to the Koszul Deformation Principle of Beilinson,
Ginzburg and Soergel (see [11], [12]).

\vskip5pt $\bullet$ Second step: the deformation is not trivial. As
mentioned before, this step can be simplified using the
C.A. Theorems.

Our proof goes exactly in the opposite direction, giving another
insight of the result: \vskip5pt $\bullet$ First step: using the
C.A. Theorems, we prove that there exists a $\CC [G]$-relative
deformation where the SRA relations hold.

\vskip5pt $\bullet$ Second step: we normalize the deformation using
an adapted equivalence defined by a symmetrization map. This is a
classical argument (e.g [8]). We then prove that the obtained
deformation is polynomial, using a powerful formula of F.A. Berezin
[6] [9].  \vskip5pt $\bullet$ Third step: we prove that the
subalgebra $(G \ * \ W) \ [\hbar]$ is isomorphic to $H_{\hbar
  \lambda}$, which is therefore a deformation, and the PBW-property.

Obviously, such a proof can be tried because we know from the
beginning what we want to find (i.e: $H_{\hbar \lambda}$), thanks to
P. Etingof and V. Ginzburg ! On the other hand, the formula of Berezin
is very explicit, and could be used to give some more light on the
structure of $H_{\hbar \lambda},$ but this is to be done.  \vskip10pt

{\bf Acknowledgments \:} \vskip5pt I am grateful to J. Alev. for
highlighting lectures on the A.F.L.S. Theorems. I thank G. Dito who
gave me an essential argument, coming from his paper [8]. I also thank
D. Arnal and R. Ushirobira for many discussions on the
subject. Finally, I am indebted to Moshe Flato for many explanations
about deformation quantization, when the theory was starting.

\vskip10pt (1) THE KOSZUL COMPLEX OF THE WEYL ALGEBRA \vskip10pt (1.1)
Let $W = \CC \ [ p_1, \ q_1,..., \ p_n, \ q_n ]$ ; there are two
algebra structures on $W$: the first one is the usual commutative
product and the second one is the Moyal $*$.product (see eg. [19] for
a short introduction); $W$, with the Moyal $*$.product, is the Weyl
Algebra. Let $V = span \ (p_1, \ q_1,..., p_n, \ q_n),$ and $\sigma
\in SP (2n)$ ; then $\sigma$ extends to an automorphism of both
algebra structures of $W.$ A Darboux basis of $V$ will be any basis of
type $\{ \sigma (p_1), \ \sigma (q_1),...,\sigma (p_n), \sigma (q_n)
\},$ for some $\sigma \in SP (2n)$.  \vskip5pt (1.2) We define
operators ${\buildrel \sim \over Z}_i$ of $W \otimes W,$ $i = 1,...,
2n,$ as follows: we set $Z_{2i-1} = p_i, \ Z_{2i} = q_i, \ i =
1...n,$ and \vskip0.5pt (1.2.1) ${\buildrel \sim \over Z}_i \ (a
\otimes b) = a * Z_i \otimes b - a \otimes Z_i \ * \ b, \ a, \ b \in
W.$

We denote by $\Lambda = \displaystyle\mathop{\oplus}_{k \geq 0}
\oalign{\cr} \Lambda^k$ the exterior algebra of $V,$ and by $i_x, \ x
\in V^*,$ the corresponding derivation of $\Lambda.$ The Koszul
complex $K = (K^k, \ d^K_k, \ k \geq - 1)$ is defined by: \vskip5pt
(1.2.2) $0 \leftarrow W \ {\buildrel {m_*} \over \longleftarrow} \ W
\otimes W \ {\buildrel {d^K_1} \over \longleftarrow}
\ ... \ {\buildrel {d^K_k} \over \longleftarrow} \ K^k = W \otimes
\Lambda^k \otimes W \leftarrow ...$ \vskip0.5pt where $m_*$ is the
Moyal product and $d^K_k := \ \displaystyle {\sum^{2n}_{i = 1}}
\ {\buildrel \sim \over Z}_i \ \otimes i_{Z^*_i}, \ k > 0.$ It is
known that $K$ is a free resolution of the bimodule $W$ (e.g:
[17]). Therefore, given a W-bimodule $M,$ applying
$\displaystyle\mathop{Hom}_{bimod} \oalign{\cr} (\bullet,M)$ on $K,$
the Hochschild cohomology $H^\bullet (M)$ is isomorphic to the
cohomology of the complex $K(M) = (K^k (M), \Delta_k,$ $k \geq 0)$:
\vskip2pt (1.2.3) $M \ {\buildrel {\Delta_\circ} \over
  \longrightarrow} \ M \otimes \Lambda^{1*} \ {\buildrel {\Delta_1}
  \over \longrightarrow} ... \rightarrow K^k (M) = M \otimes
\Lambda^{k*} = {\mathcal L} (\Lambda^k, M) \ {\buildrel {\Delta_k}
  \over \longrightarrow} ...$ with $\Delta_k = \displaystyle
       {\sum^{2n}_{i = 1}} \ ad \ Z_i \otimes \mu_{Z^*_i},$ $ad \ Z_i
       (m) = Z_i \ m - m \ Z_i, \ m \in M$ and $\mu_{Z^*_i} (\Omega) =
       Z^*_i \ \wedge \ \Omega, \ \Omega \in \Lambda^*.$ $K (M)$ will
       be called the Koszul complex of $M.$ \vskip5pt (1.3) Let
       ${\mathcal B} = ({\mathcal B}^k, \ d^{\mathcal B}, \ k \geq -
       1)$ be the normalized Bar-resolution of $W.$ Recall that
       ${\mathcal B}^{-1} = W,$ and

${\mathcal B}^k = W \otimes T^k (W / \ \CC) \otimes W,$ $k \geq 0,$
       where $T(W / \ \CC)$ is the tensor algebra of $W / \ \CC.$ We
       define an inclusion from $K^k$ into ${\mathcal B}^k$ by:

\vskip5pt (1.3.1) $a, b \in W, \ a \otimes (Z_{i_1}, \wedge...\wedge
\ Z_{i_k}) \otimes b = \displaystyle {\sum_{\sigma \in \Sigma_k}}
\ \epsilon (\sigma) \ a \otimes Z_{i_{\sigma (1)}} \otimes ... \otimes
\ Z_{i_{\sigma (k)}} \otimes b.$ The Koszul complex $K$ has a
remarkable property with respect to the normalized Bar-resolution:

\vskip5pt (1.3.2) LEMMA \vskip5pt {\it The Koszul complex $K$ is a
  subcomplex of the normalized Bar-resolution (i.e: the inclusion map
  (1.3.1) is a chain map).}  \vskip5pt $\underline {\rm Proof }$: One
checks that $d^B \vert_K = d^K$ by a straightforward direct
computation, with main argument that $[a,b] \in \CC$ if $a,b \in V.$
\hskip350pt Q.E.D.  \vskip5pt (1.3.3) {\it Remark}: a similar result
holds in the case, eg, of a polynomial algebra, and has useful
consequences. For the Weyl Algebra, it also has useful consequences,
as we shall show.  \vskip5pt (1.4) Let ${\mathcal C} (M) = ({\mathcal
  C}^k (M), \ k \geq 0, d)$ be the normalized Hochschild complex of a
bimodule $M.$ Since ${\mathcal C} (M),$ and $K(M)$ are obtained when
applying the same operation (namely $\displaystyle\mathop{Hom}_{bimod}
\oalign{\cr} (\bullet,M))$ on ${\mathcal B}$ and $K,$ by (1.3.2), one
gets a restriction map: \vskip5pt (1.4.1) $R: {\mathcal C} (M)
\longmapsto K (M),$ defined by $R (C) = C \vert_\Lambda,$ if $C \in
            {\mathcal C} (M),$ and satisfying $\Delta \circ R = R
            \circ d.$ By standard arguments ([7]), one has: \vskip10pt
            (1.4.2) PROPOSITION:

{\it The restriction map induces an isomorphism in a cohomology.}
\vskip5pt
(1.4.2) is a useful improvement of the usual isomorphism $H^\bullet (M) \simeq H^\bullet (d) \simeq H^\bullet (\Delta),$ since the 
isomorphism is explicit: it comes from the restriction map. For instance, one has the following immediate consequence:
\vskip5pt
(1.4.3) COROLLARY

{\it Let $c \in {\mathcal L} (\Lambda^k, M)$ be a Koszul cocycle of
  $M,$ then there exists a Hoschild cocycle $C \in {\mathcal C}^k (M)$
  such that $C \vert_\Lambda^k = c$. If a second Hochschild cocycle
  $C'$ has the same property, then $C - C'$ is a Hochschild
  coboundary.}

\vskip5pt
(1.5) {\it Remark}: we have defined the Koszul resolution using the canonical Darboux basis ${p_1, q_1,..., p_n, q_n },$ let us show that the Koszul resolution has an intrinsic nature, so that the formulas are valid in any basis of $V.$ To do that, we define, for any $X \in V,$ an operator $\tilde{X}$ of $W \otimes W$ by $\tilde{X} (a \otimes b) = a * X \otimes b - a \otimes X * b, \ \forall a,b \in W.$ Then we define $\rho_k: {\mathcal L} (V) \mapsto {\mathcal L} (K^k, \ K^{k-1})$ by $\rho_k (X \otimes \varphi) = \tilde{X} \otimes i_\varphi, \ X \in V, \ \varphi \in V^*.$ Since $d^K_k = \rho_k (Id_V),$ the result follows. We shall use this remark in section 2.
\vskip10pt
(2) HOCHSCHILD COHOMOLOGY OF THE TWISTED BIMODULE $W_\sigma$
\vskip10pt
(2.1) Given any automorphism $\sigma$ of the Weyl Algebra, we denote by $W_\sigma$ the W-bimodule with underlying space $W$ and action:
\vskip5pt
(2.1.1) $a,b \in W, a \ \displaystyle\mathop{\bullet}_\sigma \oalign{\cr} \ b = a * b,$ $b \ \displaystyle\mathop{\bullet}_\sigma \oalign{\cr} \ a = b * \sigma (a).$

Let ${\mathcal C} (W_\sigma)$ be the normalized Hochschild complex of $W_\sigma$ ; then ${\mathcal C}^k (W_\sigma) = C^k (W) \ = \ {\mathcal L} (\displaystyle\mathop{\otimes}_k \oalign{\cr} \ W / \CC, W).$ Let $H^\bullet (W_\sigma)$ be the corresponding Hochschild cohomology.
\vskip5pt

(2.1.1) When $\sigma = Id_W,$ one obtains the usual Hochschild
cohomology of $W,$ which is well known ([20], [10]): $H^0 (W) = \CC,$
and $H^k (W) = {0},$ if $k > 0.$ As a consequence, $W$ is rigid in
Gerstenhaber deformation theory.

\vskip5pt (2.1.2) When $\sigma$ is the parity of $W$, $H^\bullet
(W_\sigma)$ was computed in [17] ; one has $dim \ H^{2n} (W_\sigma) =
1,$ and $H^k (W_\sigma) = \{ 0 \},$ if $k \neq 2n.$ As a consequence,
if $n > 1,$ $W$ is rigid in super-commutative deformation theory (see
[16] [17]) ; when $n=1$ the enveloping algebra ${\mathcal
  U}(osp(1,2))$ provides a universal super-commutative deformation,
with deformation parameter the ghost of ${\mathcal U}(osp(1,2))$
([17], [18], [4]).  

\vskip10pt (2.2) Let us assume that $\sigma$ is the automorphism of
$W$ extending a diagonalizable element $\sigma$ of $SP (2n).$ We
assume moreover that $\sigma \not = Id,$ and introduce $V_\sigma =
{\rm range} \ (\sigma - Id) \vert_V.$ Then $dim \ V_\sigma$ is even,
say $2 k_\sigma.$ (see [2])

Let ${x_1,...,x_{2k_{\sigma}}}$ be a Darboux basis of $V_\sigma,$ and
$\omega_\sigma \:= \ x^*_1 \ \wedge ... \wedge \ x^*_{2 k_\sigma}$
(any choice of the Darboux basis will lead to the same
$\omega_\sigma$) ; if $\omega$ is the canonical two form, one has:

\vskip5pt (2.2.1) $\omega_\sigma (v_1,...,v_{2 k \sigma}) =
\frac{1}{k_\sigma !} \omega^{k \sigma} (P_\sigma (v_1),...,P_\sigma
(v_{2 k \sigma})),$ where $P_\sigma$ is the projection on $V_\sigma$
associated to $V = V_\sigma \oplus Ker (\sigma - Id)$ (see [2] [3]
[11] for details). One has: \vskip10pt (2.2.2) C.A. THEOREM 1
\vskip5pt {\it (1) [1] $dim H^{2k_\sigma} (W_\sigma) = 1,$ and $H^k
  (W_\sigma) = \{ 0 \},$ if $k \not = 2 k_\sigma.$ \vskip5pt (2) There
  exists a Hochschild cocycle $\Omega_\sigma \in Z^{2k_\sigma}
  (W_\sigma)$ such that:} $\Omega_\sigma \vert_{\Lambda^{2k_\sigma}} =
\omega_\sigma, \ {\rm and} \ H^{2k_\sigma} (W_\sigma) \ = \ \CC
. \Omega_\sigma.$ \vskip5pt {\it (3) Let $S$ be a finite subgroup of
  $SP (2n),$ commuting with $\sigma.$ Denoting by $H^\bullet_S
  (W_\sigma)$ the $S$-invariant Hochschild cohomology (i.e: computed
  from $S$-invariant cochains), one has $H^\bullet_S (W_\sigma) =
  H^\bullet (W_\sigma),$ and the cocycle of (2) can be chosen
  $S$-invariant.}  \vskip10pt $\underline {\rm Proof }$: \vskip5pt By
      [2], there exists a Darboux basis ${P_1, Q_1,...,P_{n}, Q_{n}}$
      of $V$ such that: \vskip0.5pt $\ \ \ \sigma (P_i) = \alpha_i
      . P_i, \ \sigma (Q_i) = \alpha^{-1}_i . Q_i, \ {\rm with}
      \ \alpha_i \not = 1, \ {\rm if } \ i \leq k_\sigma, \ {\rm and}
      \ \alpha_i = 1, \ {\rm if} \ i > k_\sigma.$ 

\vskip5pt We compute $H^\bullet (W_\sigma)$ using the Koszul complex
(see remark (1.5)). Using (1.2.3), and Moyal product, the differential
is given by:

\vskip5pt (2.2.3) $\Delta_\sigma = \displaystyle {\sum^{2n}_{i = 1}}
\ T_i \otimes \mu_{Z^*_i},$ where $T_{2i-1} = (1 - \alpha_i) m_{P_i} +
   {1 \over 2} (1 + \alpha_i) \ {\partial \over \partial Q_i},$
   $T_{2i} = (1 - \alpha^{-1}_i) m_{Q_i} - {1 \over 2} (1 +
   \alpha^{-1}_i) {\partial \over \partial P_i},$ $m_{P_i} (a) = P_i
   \ a, \ m_{Q_i} (a) = Q_i\ a, \ a \in W.$

\vskip5pt We define operators $\theta$ and $A$ of $W$ by:
\[ \ \ \ \ \ (2.2.4) \begin{cases} A (P_i) = (1 - \alpha_i)^{-1} \ . \ P_i,& \\ &{\rm if} \ i \leq k_\sigma, \\ A (Q_i) = (1 - \alpha^{-1}_i)^{-1} \ . \ Q_i & \end{cases} \quad \text{ and } \begin{cases} A (P_i) = Q_i,& \\ &{\rm if} \ i > k_\sigma, \\ A (Q_i) = -P_i & \end{cases}\]

\vskip0.1pt (2.2.5) $\ \ \ \theta = exp [ - {1 \over 2}
  \ \displaystyle {\sum_{i \leq k_\sigma}} \ \beta_i \ {{\partial^2}
    \over {\partial \ P_i . \partial \ Q_i}}],$ where $\beta_i =
       {{1+\alpha_i} \over {1 - \alpha_i}}.$ \vskip0.5pt Let now $\xi
       \:= \ A \circ \theta,$ and $\Delta' \:= \ (\xi \otimes
       Id_{\wedge^*}) \ \circ \ \Delta_\sigma \ \circ \ (\xi \otimes
       Id_{\wedge^*})^{-1}.$ It is easy to check that: \vskip5pt
       (2.2.6) $\Delta' = \displaystyle {\sum^{2k_\sigma}_{i = 1}}
       \ m_{Z_i} \otimes \mu_{Z^*_i} + \displaystyle {\sum^{2n}_{i = 2
           k_\sigma + 1}} \ {{\partial} \over {\partial Z_i}} \otimes
       \mu_{Z^*_i}.$

Since $\Delta^{'2} = 0,$ we get a new complex, with cohomology
$H^\bullet (\Delta'),$ and from the definition of $\Delta'$, $\xi^{-1}
\otimes Id: Z(\Delta') \rightarrow Z(\Delta_\sigma)$ induces an
isomorphism $H (\Delta') \simeq H (\Delta_\sigma).$ So, we get a new
equivalent differential $\Delta' = \Delta'_1 + \Delta'_2,$ where
$\Delta'_1$ corresponds to the case $\alpha_i = - 1,$ $\forall i,$ and
$\Delta'_2$ is a de-Rham type differential.  \vskip5pt (2.2.7) It is
obvious that $\omega_\sigma \in Z^{2 k_\sigma} (\Delta')$ ; given $C
\in B^{2 k_\sigma} (\Delta'),$ one has \linebreak $C
(Z_1,...,Z_{2k_\sigma}) (0) = 0,$ and since $\omega_\sigma
(Z_1,...,Z_{2k_\sigma}) = 1,$ $\omega_\sigma \notin B^{2 k_\sigma}
(\Delta').$ Moreover $(\xi^{-1} \otimes Id) (\omega_\sigma) =
\omega_\sigma,$ so $\omega_\sigma$ is also a non trivial cocycle in
$Z^{2 k_\sigma} (\Delta_\sigma).$ By (1.4.3), there exists a non
trivial Hochschild cocycle $\Omega_\alpha \in Z^{2k_\sigma}
(W_\sigma)$ such that $\Omega_\sigma \ \vert_{\Lambda^{2k_\sigma}} =
\omega_\sigma.$ \vskip5pt (2.2.8) We need some more notations. Let
$V_1 = V_\sigma,$ $V_2 = Ker (\sigma - Id),$ $\Lambda_1$ and
$\Lambda_2$ the exterior algebras of $V_1$ and $V_2$, $W_1$ and $W_2$
the symmetric algebras of $V_1$ and $V_2.$ One has $\Lambda =
\Lambda_1 \otimes \Lambda_2,$ and $W = W_1 \otimes W_2.$ We introduce
the following operators: \vskip2pt $h_1 = \displaystyle
{\sum^{2k_\sigma}_{i = 1}} \ {\partial \over \partial Z_i} \otimes
i_{Z_i}, \ h_2 = \displaystyle {\sum_{i > 2k_\sigma}} \ m_{Z_i}
\otimes i_{Z_i},$ $R_1 = \displaystyle {\sum_{i \leq 2k_\sigma}} \ Z_i
\ {\partial \over \partial Z_i}, \ R_2 = \displaystyle {\sum_{i >
    2k_\sigma}} \ Z_i \ {\partial \over \partial Z_i},$ ${\mathcal
  R}_1 = \displaystyle {\sum_{i \leq 2k_\sigma}} \ Z^*_i \ \wedge
i_{Z_i},$ $\ \ \ \ \ \ {\mathcal R}_2 = \displaystyle {\sum_{i >
    2k_\sigma}} \ Z^*_i \ \wedge \ i_{Z_i}.$ \vskip5pt (2.2.9) The
complex $(K(W_\sigma), \Delta')$ splits into three sub-complexes:
\vskip5pt $K(W_\sigma) = \Lambda^{2 k_\sigma *}_1 \oplus H_1 \oplus
H_2 \ ,$ $H_1 = W_1 \otimes (\displaystyle\mathop{\oplus}_{i <
  2k_\sigma} \ \Lambda^{i*}_1) \oalign{\cr} ) + W^{+}_1 \otimes
\Lambda^{2 k_\sigma *}_1,$ $H_2 = W_1 \otimes \Lambda^*_1 \otimes
\Lambda^{+*}_2 + W_1 \otimes W^+_2 \otimes \Lambda^*_1 \otimes
\Lambda^*_2,$ where the subscript + stands for the kernel of the
corresponding evaluation ; (it is easy to see that $H_1$ and $H_2$ are
$\Delta'$-stable, and $\bigwedge^{2 k_\sigma *}$ is a subcomplex by
(2.2.7)).  \vskip5pt (2.2.10) By a simple computation, one has $h_2
\circ \Delta' + \Delta' \circ h_2 = R_2 + {\mathcal R}_2 \:= \ T_2.$
Since $\tau_2 \:= \ T_2 \ \vert_{H^2}$ is invertible, since moreover
$\Delta' \circ \tau_2 = \Delta' \circ h_2 \circ \Delta' = \tau_2 \circ
\Delta',$ we obtain that $Id_{H_2} = (h_{2^\circ} \circ \tau^{-1}_2)
\circ \Delta' + \Delta' \circ (h_2 \circ \tau^{-1}_2),$ so the complex
$(H_2, \Delta')$ has trivial cohomology.  \vskip5pt (2.2.11) Let
$\Delta'_1 = \displaystyle {\sum_{i \leq 2k_\sigma}} \ m_{Z_i} \otimes
\mu_{Z^*_i} = \Delta' \ \vert_{H_1}.$ A computation gives $h_1 \circ
\Delta'_1 + \Delta'_1 \circ h_1 = R_1 + 2k_\sigma . Id - {\mathcal
  R}_1 \:= \ T_1.$ Let $\tau_1 = T_1 \ \vert_{H_1},$ $\tau_1$ commutes
with $\Delta'_1,$ so $(h \circ \ \tau^{-1}_1) \ \circ \Delta'_1 +
\Delta'_1 \ \circ (h_1 \circ \tau^{-1}_1) = Id_{H_1}$ ; it results
that the complex $(H_1, \Delta')$ has trivial cohomology.  \vskip5pt
(2.2.10) Finally, we stay with the complex $(\Lambda^{2 k_{\sigma} *},
\ \Delta' = 0)$ treated in (2.2.7): $\omega_\sigma$ is a non trivial
cocycle; so finally $H^\bullet (\Delta') = \CC.\omega_\sigma.$
Applying (2.2.6) and (2.2.7), it follows that $H^\bullet
(\Delta_\sigma) = \CC . \omega_\sigma.$ \vskip5pt (2.2.11) By (1.4.2)
and (1.4.3), there exists a Hochschild cocycle $\Omega_\sigma \in Z^{2
  k_\sigma} \ (W_\sigma)$ such that the Hochschild cohomology
$H^\bullet (W_\sigma) = \CC . \Omega_\sigma,$ and $\Omega_\sigma
\ \vert_{\Lambda^{2 k_\sigma}} = \omega_\sigma.$ \vskip5pt (2.2.11)
Let $S$ be a finite subgroup of $SP (2n)$ such that $s \ \circ \sigma
= \sigma \ \circ s, \ \forall \ s \in S.$ $S$ acts on Hochschild
cochains $C \in {\mathcal C}^k (W)$ by:

$\pi_s (C) (a_1,...,a_k) = s (C (s^{-1} (a_1),..., s^{-1} \ (a_k))),
\ a_i \in W,$ and this action commutes with the Hochschild
differential. So we can consider the $S$-invariant Hochschild
cohomology $H^\bullet_S (W_\sigma),$ using the complex ${\mathcal
  C}^\bullet_S (W_\sigma),$ of $S$-invariant cochains. By standard
arguments, since $S$ is finite, one has an inclusion $H^\bullet_S
(W_\sigma) \subset H^\bullet (W_\sigma).$ We now show that this is an
equality:

\vskip5pt (2.2.12) Let $P = {1 \over \vert S \vert} \ \displaystyle
          {\sum_{s \in S}} \ s,$ then $\pi_P$ is a projector from
          ${\mathcal C} (W)$ onto its trivial isotopic component
          \hskip5pt ${\mathcal C}_S (W).$ So $\pi_P (\Omega_\sigma)$
          is an $S$-invariant Hochschild cocycle, and one has $\pi_P
          (\Omega_\sigma) \ \vert_{\Lambda^{2 k_\sigma}} = \pi_P
          (\omega_\sigma) = \omega_\sigma.$ (because $\pi_s
          (\omega_\sigma) = \omega_\sigma, \ \forall \ s \in S,$ since
          $s \in SP (V_\sigma)).$ Using (1.4.2), $H^{2 k_\sigma}
          (W_\sigma) = \CC . \pi_P (\Omega_\sigma) = H^{2 k_\sigma}_S
          (W_\sigma),$ and $H^k_S (W_\sigma) = H^k (W_\sigma) = \{ 0
          \}$ if $k \not = 2 k_\sigma.$ Therefore $H^\bullet
          (W_\sigma) = H^\bullet_S (W_\sigma) = \CC . \pi_P
          (\Omega_\sigma).$ \hfill Q.E.D.

\vskip5pt (2.3) {\it Remark }: Let $\sigma, \ V_\sigma, \ k_\sigma,
\ \omega_\sigma, \ \Omega_\sigma$ be as in (2.2), and $\tau = x
\ \sigma \ x^{-1},$ $x \in SP (2n)$ ; we introduce corresponding
$V_\tau = \pi_x (V_\sigma), \ k_\tau = k_\sigma, \ \omega_\tau = \pi_x
(\omega_\sigma).$ Let $d_\sigma$ and $d_\tau$ be the respective
Hochschild differentials of the bimodules $W_\sigma$ and $W_\tau$,
since $d_\tau = \pi_x \ \circ \ d_\sigma \ \circ \ \pi_{x^{-1}},$ one
has, for Hochschild cocycles, $Z (W_\tau) = \pi_x (Z (W_\sigma)),$ the
same for coboundaries, and $\pi_x$ induces an isomorphism from
$H^\bullet (W_\sigma)$ onto $H^\bullet (W_\tau)$. The cocycle $\pi_x
(\Omega_\sigma)$ is a non trivial element in $Z (W_\tau)$, which
satisfies $\pi_x (\Omega_\sigma) \ \vert_{\Lambda^{2 k_\sigma}} =
\omega_\tau,$ so, applying (2.2.2) (2) to $\tau,$ we can choose
$\Omega_\tau = \pi_x (\Omega_\sigma).$ Let $S$ be a finite subgroup of
$SP (2n)$, if $S$ commutes with $\sigma,$ then $x S x^{-1}:= \ S'$
commutes with $\tau,$ and a cochain $C$ is $S$-invariant if and only
if $\pi_x (C)$ is $S'$-invariant. It results that $H^\bullet_{S'}
(W_\tau) = \pi_x (H^\bullet_S (W_\sigma)) = \pi_x (H^\bullet
(W_\sigma)) = H^\bullet (W_\tau)$ ; moreover, using (2.2.2.)  (3), we
can start with an $S$-invariant $\Omega_\sigma,$ then $\Omega_\tau =
\pi_x (\Omega_\sigma)$ is $S'$-invariant, and $H^\bullet (W_\tau) =
\CC . \Omega_\tau.$ \vskip10pt (2.4) {\it Remark }: We develop (2.3)
in a context which will be used in the next section: Let $G$ be a
finite subgroup of $SP (2n).$ Any $g \in G$ satisfies $g^{card \ G} =
1,$ so $g$ is diagonalizable, and we can apply all the results of
(2.2) and (2.3). Given a conjugacy class $\gamma$ of $G,$ we fix
$\sigma \in \gamma,$ and use the notations of (2.2) and (2.3). If
$\tau \in \gamma,$ one has $k_\tau = k_\sigma$ (denoted by $k_\gamma$)
and $\omega_\tau = \pi_x (\omega_\sigma)$ for any $x$ such that $\tau
= x \ \sigma \ x^{-1}$ (2.3). Denote by $S_\tau$ the centralizer of
$\tau \in \gamma$ in $G.$ One has $S_\tau = x \ S_\sigma \ x^{-1},$ if
$\tau = x \ \sigma \ x^{-1},$ and by (2.2.2) $H^\bullet_{S_\tau}
(W_\tau) = H^\bullet (W_\tau),$ the cohomology being one-dimensional,
concentrated in degree $2 \ k_\gamma.$ Starting with an
$S_\sigma$-invariant $\Omega_\sigma,$ given by (2.2.2) (3), we define
$\Omega_\tau = \pi_x (\Omega_\sigma),$ if $\tau = x \ \sigma \ x^{-1}
\in \gamma$ ; first, we remark that if $\tau = x' \ \sigma \ x'^{-1},$
then $\pi_{x'} (\Omega_\sigma) = \pi_x (\Omega_\sigma),$ secondly, by
(2.3), $\Omega_\tau$ is $S_\tau$-invariant, $H^{2 k_\tau} (W_\tau) =
\CC . \Omega_\tau,$ and one has $\Omega_\tau \ \vert_{\Lambda^{2
    k_\gamma}} = \omega_\tau.$ \vskip10pt (3) COHOMOLOGY OF $G \ *
\ W$ \vskip10pt (3.1) In this section, $G$ is a finite subgroup of $SP
(2n)$, $G \ * \ W$ is the algebra with underlying space $W \ \otimes
\ \CC [G],$ and product such that: \vskip5pt (3.1.1) $\sigma \ * \ a =
\sigma (a) \otimes \sigma, \ a \ * \ \sigma = a \otimes \sigma,
\ \forall \ \sigma \in G, \ a \in W, \ {\rm Moyal \ product \ on} \ W
= W \otimes 1,$ and convolution product on $\CC [G] = 1 \otimes \CC
[G].$ The group $G$ acts on $G * W$ by conjugation: \vskip5pt (3.1.2)
$Ad \sigma (b) = b * \sigma * b^{-1}, \ \forall \sigma \in G, \ b \in
G * W,$ and the product is invariant: \vskip5pt (3.1.3) $Ad \sigma
(b_1 * b_2) = Ad \sigma (b_1) * Ad \sigma (b_2).$ On $W,$ one has: $Ad
\sigma (a) = \sigma (a),$ and $Ad \sigma (a_1 * a_2) = \sigma (a_1) *
\sigma (a_2), \ \forall \ \sigma \in G, \ a, a_1, a_2 \in W.$
\vskip5pt (3.2) $\CC [G]$ is a separable algebra, so $H^\bullet (G *
W) = H^\bullet_{\CC [G]} (G * W),$ [13], hence Hochschild cohomology
of $G * W$ can be computed using $\CC [G]$-relative cochains, i.e
cochains $C$ which satisfy:

\vskip5pt (3.2.2) $C (b_1,...,b_k) = 0,$ if some $b_i \in \CC [G],
\ \sigma * C (b_1,...,b_k) = C (\sigma * b_1, b_2,...,b_k),$ $C
(b_1,...,b_i * \sigma, b_{i+1},...,b_k)= C (b_1,...,b_i, \ \sigma *
b_{i+1},...,b_k), \ \forall \ i,$ and $C (b_1,...,b_k * \sigma) = C
(b_1,...,b_k) \ * \ \sigma, \ \forall \ \sigma \in G, \ b_i \in G *
W.$

It results that: \vskip5pt (3.2.3) $Ad \sigma (C(b_1,...,b_k)) = C(Ad
\sigma (b_1),..., Ad \sigma (b_k)), \ \forall \ \sigma \in G, \ b_i
\in G * W,$ and, as a particular case, that if $D \:= \ C \ \vert
\otimes_k W,$ one has:

\vskip5pt
(3.2.4) $Ad \sigma (D(a_1,...,a_k)) = D (\sigma(a_1),..., \sigma (a_k)), \ \forall \ \sigma \in G, \ a_i \in W,$ i.e: $D$ is $G$-invariant.

Conversely, given $D \in {\mathcal L}(\otimes_k W, G *W),$ $G$
invariant and normalized, one defines a \ $\CC [G]$-relative cochain
$C$ of $G * W,$ such that $C \vert \otimes_k W = D,$ by: \vskip5pt
(3.2.5) $C (a_1 \otimes g_1,..., a_k \otimes g_k) \:= \ D(a_1, g_1
(a_2),..., g_1...g_{k-1} (a_k)) * g_1...g_k.$

So we can replace $\CC [G]$-relative $k$-cochains on $G * W$ by their
restrictions to $\otimes_k \ W,$ and this is what we shall now do. We
denote by ${\mathcal C}^k_{ \CC [G]} \ (G * W)$ the corresponding
space, i.e $C \in {\mathcal L} (\ \otimes_k \ W, \ G *W)$ such that:
\vskip5pt (3.2.6) $C$ is normalized: $C (a_1,...,a_k) = 0,$ if some
$a_i \in \ \CC,$ and $C$ is $G$-invariant.

Let us quickly recall what is a $\CC [G]$-relative deformation of the product of $G * W.$ Denote by $m_*$ the product of $G * W$, and $m_\hbar = m_* + \displaystyle {\sum_{k \geq 1}} \ \hbar^k . m_k$ the new-product. $m_\hbar$ is $\CC [G]$-relative if all cochains $m_k$ are $\CC [G]$-relative cochains. Observe that the product $m_*$ is $G$-invariant by (3.1.3). It results that $m_\hbar$ is $\ \CC [G]$-relative if and only if 
\vskip5pt
(3.2.7) $Ad \sigma (m_\hbar(a_1, a_2)) = m_\hbar(\sigma (a_1), \sigma (a_2)), \ \forall \ \sigma \in G, a_1, a_2 \in W$ ; in that case, one has:
\vskip5pt
(3.2.8) $Ad \sigma (m_\hbar(b_1, b_2)) = m_\hbar(Ad \sigma (b_1), Ad \sigma (b_2)), \ \forall \ \sigma \in G, b_1, b_2 \in G * W$ ; i.e $m_\hbar$ is $G$-invariant.

\vskip5pt (3.3) When $C \in {\mathcal C}^k_{ \CC [G]} (G * W)$ we
write $C = \displaystyle {\sum_{g \in G}} \ C_g \otimes g,$ with $C_g
\in {\mathcal C}^k (W),$ the space of normalized cochains on $W.$ We
recall that $G$ acts on ${\mathcal C}^k (W)$ by:

\vskip5pt
(3.3.1) $\pi_\sigma (C) (a_1,...,a_k) = \sigma (C(\sigma^{-1} (a_1),...,\sigma^{-1} (a_k))).$ \hskip10pt The $G$-invariance condition becomes:
\vskip5pt
(3.3.2) $C_{Ad \sigma (g)} = \pi_\sigma (C_g), \ \forall \ \sigma, g \in G.$

Let $S_\sigma$ be the centralizer of $\sigma$ in $G$, observe that $\pi_s (C_\sigma) = C_\sigma, \ \forall s \in S_\sigma,$ so that $C_\sigma$ has to be $S_\sigma$-invariant. We shall use the notation:
\vskip5pt
(3.3.3) ${\mathcal C}^k_\sigma (W) \:= \ {\mathcal C}^k (W)^{S_\sigma}.$

Let $\Gamma$ be the set of conjugacy classes of $G$ ; we fix a section $\sigma_\gamma \ \in \ \gamma, \ \forall \ \gamma \ \in \ \Gamma,$ and, with an abuse of notation, we write $\gamma = \{ \bar{x} \in G/S_\gamma \},$ $S_\gamma$ being the centralizer of $\sigma_\gamma$ in $G.$ Given $C \in {\mathcal C}^k_{ \CC [G]} \ (G * W),$ let $\tilde{C}_\gamma:= C_{\sigma_\gamma},$ $\gamma \in \Gamma,$ then, from (3.3.2), $\tilde{C}_\gamma \in {\mathcal C}^k_\gamma (W):= {\mathcal C}^k_{\sigma_\gamma} (W),$ and, using (3.3.2), one gets:
\vskip5pt
(3.3.4) $C = \displaystyle {\sum_{\gamma \in \Gamma}} \ \displaystyle {\sum_{\bar{x} \in G/S_\gamma}} \ \pi_x (\tilde{C}_\gamma) \otimes Ad \ x (\sigma_\gamma).$

On the other hand, given $\tilde{C} = (\tilde{C}_\gamma, \ \gamma \in
\Gamma) \in \ \displaystyle\mathop{\Pi}_{\gamma \in \Gamma}
\oalign{\cr} \ {\mathcal C}^k_\gamma (W),$ the map $B(x):= \pi_x
(\tilde{C}_\gamma), \ x \in G,$ satisfies $B(x s) = B(x), \ \forall
\ s \in S_\gamma,$ so, if we define $C$ by formula (3.3.4), then $C
\in {\mathcal C}^k_{ \CC [G]} (G * W).$ So we have proved:

\vskip10pt {(3.3.5) LEMMA \vskip5pt {\it The map $T: C \rightarrow
    \tilde{C}$ is an isomorphism from ${\mathcal C}^k_{ \CC [G]} (G *
    W)$ onto $\displaystyle\mathop{\Pi}_{\gamma \in \Gamma}
    \oalign{\cr} \ {\mathcal C}^k_\gamma (W).$} \vskip5pt (3.4) Let
  $d$ and $d_g, \ g \in G,$ be the differentials of the Hochschild
  complex respectively of $G * W$ and $W_g.$ One has: \vskip5pt
  (3.4.1) If $C = \displaystyle {\sum_{g \in G}} \ C_g \otimes g \in
  {\mathcal C}^k_{ \CC [G]} (W * G),$ then $d(C) = \displaystyle
  {\sum_{g \in G}} \ d_g (C_g) \otimes g.$ \vskip5pt {(3.4.2) LEMMA
    \vskip5pt {\it With the notations of (3.4.1) and (3.3.4), $C \in
      B^k_{ \CC [G]} (G * W)$ if and only if $\tilde{C}_\gamma \in
      B^k_\gamma (W_\gamma), \ \forall \gamma,$ where $W_\gamma:=
      W_{\sigma_\gamma}, \ \gamma \in \Gamma.$} \vskip5pt $\underline
              {\rm Proof }$: Let us assume that $C \in Z^k_{ \CC [G]}
              (G * W)$, and that $\tilde{C}_\gamma \in B^k_\gamma
              (W_\gamma), \ \forall \ \gamma \in \Gamma.$ One has $C =
              \displaystyle {\sum_{\gamma \in \Gamma}} \ \displaystyle
                            {\sum_{\bar{x} \in G/S_\gamma}} \ \pi_x
                            (\tilde{C}_\gamma) \otimes Ad \ x
                            (\sigma_\gamma).$ $\tilde{C}_\gamma$ is
                            $S_\gamma$-invariant, therefore, by
                            (2.2.2) there exists an
                            $S_\gamma$-invariant cochain
                            $\tilde{B}_\gamma$ such that
                            $\tilde{C}_\gamma = d_{\sigma_\gamma}
                            \tilde{B}_\gamma$. Using (2.3), one has
                            $\pi_x (\tilde{C}_\gamma) = d_{Ad x
                              (\sigma_\gamma)} \ [\pi_x
                              (\tilde{B}_\gamma)], \ \forall \ x \in
                            G$ ; defining $B := \displaystyle
                            {\sum_{\gamma \in \Gamma}} \ \displaystyle
                            {\sum_{\bar{x} \in G/S_\gamma}} \ \pi_x
                            (\tilde{B}_\gamma) \otimes Ad \ x
                            (\sigma_\gamma),$ $B$ is a $\CC
                            [G]$-relative cochain by (3.3.5), and one
                            has $C = d B$ by (3.4.1) \ Q.E.D.

\vskip10pt (3.4.3) PROPOSITION \vskip5pt {\it (1) The map $T$ of
  (3.3.5) induces an isomorphism from} $H^k_{ \CC [G]} \ (G * W)$ onto
\linebreak $\displaystyle\mathop{\Pi}_{\gamma \in \Gamma} \oalign{\cr}
\ H^k_\gamma (W_\gamma) = \ \displaystyle\mathop{\Pi}_{\gamma \in
  \Gamma} \oalign{\cr} \ H^k \ (W_\gamma)$

(2) [1] $H^k (G * W) = \{ 0 \},$ if $k$ is odd. Let $\Gamma_{2k} = { \gamma \in \Gamma \ / \ k_\sigma = k, \ \forall \ \sigma \in \gamma},$ then $dim \ H^{2k} (G * W) = card \ \Gamma_{2k}.$}
\vskip5pt
$\underline {\rm Proof}$
\vskip5pt
(1) We need only to prove that $T$ is onto, but this is an immediate consequence of formula: $d_{Ad x (\sigma_\gamma)}$ $(\pi_x (\tilde{C}_\gamma)) = \pi_x \ (d_{\sigma_\gamma} (\tilde{C}_\gamma)),$ which is proved in (2.3), so if $C$ is defined by (3.3.4), $d \ C = 0.$

(2) We apply (1) and (2.2.2). \hskip100pt Q.E.D.  

\vskip5pt (3.4.4) {\it Remark}: As proved in [1], by dimension
argument, there is an isomorphism $H^{2k} (G * W) \simeq
\CC^{\Gamma_{2k}}.$ \vskip5pt (3.5) Let us precise the isomorphism of
(3.4.4). We assume that $\Gamma_{2k} \neq \phi,$ and take $\gamma \in
\Gamma_{2k}.$ Using (2.2.2), there exists $\Omega_\gamma \in
Z^{2k}_{S_\gamma} (W_\gamma)$ such that $H^{2k} (W_\gamma) = \CC
. \Omega_\gamma$ and $\Omega_\gamma \ \vert_{\Lambda^{2k}} =
\omega_{\sigma_\gamma}.$ It results from (2.3) that for any $\tau = Ad
x (\sigma_\gamma) \in \gamma,$ one has $\pi_x (\Omega_\gamma)
\vert_{\Lambda^{2k}} = \omega_\tau,$ and $H^{2k} (W_\tau) = \CC
. \pi_x (\Omega_\gamma).$ Let us define: ${\bf C}_\gamma:=
\displaystyle {\sum_{\bar{x} \in G/S_\gamma}} \ \pi_x (\Omega_\gamma)
\otimes Ad x (\sigma_\gamma),$ and decompose: ${\bf C}_\gamma =
\displaystyle {\sum_{g \in G}} \ C^\gamma_g \otimes g,$ then one has
$C^\gamma_g = 0,$ if $g \notin \gamma,$ $C^\gamma_g
\vert_{\Lambda^{2k}} = \omega_g,$ if $g \in \gamma,$ and $\tilde{\bf
  C}_\gamma = C^\gamma_{\sigma_\gamma} = \Omega_\gamma.$ \vskip5pt
(3.5.1) PROPOSITION \vskip5pt {\it (1) ${{\bf C}_\gamma, \ \gamma \in
    \Gamma_{2k}}$ is a basis of $H^{2k}_{ \CC [G]} (G * W).$}

{\it (2) Given $\lambda \in \CC^{\Gamma_{2k}},$ there exists a cocycle $C_\lambda \in Z^{2k}_{ \CC [G]} (G * W)$ such that:}

$\ \ \ C_\lambda (X_1 \wedge ... \wedge X_{2k}) = \displaystyle {\sum_{\gamma \in \Gamma_{2k}}} \ \lambda (\gamma) \ \displaystyle {\sum_{g \in \gamma}} \ \omega_g (X_1,..., X_{2k}) \otimes g, \ \forall X_i \in V, \ g \in G.$

$C_\lambda$ {\it is a coboundary if and only if} $\lambda = 0,$ {\it and the map} $\lambda \rightarrow C_\lambda$ {\it induces an isomorphism from} $\CC^{\Gamma_{2k}}$ {\it onto} $H^{2k}_{\CC [G]} (G * W).$}

{\it (3) If a cocycle $C \in Z^{2k}_{\CC [G]} (G * W)$ vanishes on
  $\Lambda^{2k},$ then $C$ is a coboundary.}  \vskip5pt $\underline
{\rm Proof }$: {\rm To obtain (1), we apply (3.4.3) ; then, we define
  $C_\lambda = \displaystyle {\sum_{\gamma \in \Gamma_{2k}}} \ \lambda
  (\gamma) \ {\bf C}_\gamma,$ and prove (2) using (3.5) and (3.5.1)
  (1). (3) is consequence of (3.4.2) and (1.4.3). \hskip15pt Q.E.D.}
\vskip5pt (3.5.2) {\it Remark}: (3.4.3) and (3.5.1) give the
C.A. Theorem 2.  \vskip10pt (4) AN ALTERNATIVE PROOF OF A THEOREM OF
P. ETINGOF AND V. GINZBURG ABOUT SYMPLECTIC REFLECTION ALGEBRAS

\vskip10pt (4.1) With the notations of (2.2), a symplectic reflection
(in the finite subgroup $G$ of $SP(2n)$), is an element $g$ such that
$dim \ V_g = 2.$ When there are no symplectic reflections in $G,$ by
(3.4.3), $G * W$ is rigid. So let us assume the contrary. Given any
$\lambda \in \CC^{\Gamma_2}, \ \lambda \neq 0,$ we construct a
non-trivial $2$-cocycle $C_\lambda$ by (3.5.1) (2) ; since $H^3_{ \CC
  [G]} (G * W) = \{ 0 \},$ $C_\lambda$ is not obstructed, so there
exists a $\CC [G]$-relative non trivial deformation of $G * W,$ with
leading cocycle $C_\lambda$. Using once more $H^3_{\CC [G]} (G * W) =
\{ 0 \},$ and varying $\lambda,$ the procedure will provide a
universal deformation formula of $G * W$ (see [11]). Using (3.5.1),
one has:

$\ \ \ \ \forall \ X,Y \in V, \ [X,Y]^\lambda_\hbar = \omega (X,Y) + \hbar \ \displaystyle {\sum_{\gamma \in \Gamma_2}} \ \lambda (\gamma) \ \displaystyle {\sum_{g \in \gamma}} \ \omega_g (X,Y) \otimes g + \hbar^2 (...).$

At first order, we find exactly the relations of the Symplectic
Reflection Algebra $H_{\hbar \lambda}$ of [12]. We call these
relations the SRA-relations.  \vskip5pt (4.1.2) By construction, since
$C_\lambda$ is non trivial, the corresponding deformation is non
trivial.  \vskip10pt (4.2) Let us quickly recall the definition of the
symplectic reflection algebra $H_{\hbar \lambda},$ following [11]
[12]. Let $T(V)$ be the tensor algebra of $V,$ $G * T(V)$ the smash
product of $T(V)$ with $\CC [G],$ and, given $\lambda \in
\CC^{\Gamma_2},$ $I_\lambda$ the ideal in [$G * T(V)$] [$\hbar$]
generated by:

$R_{\hbar \lambda} (X,Y) = X \otimes Y - Y \otimes X - \omega (X,Y) - \hbar \ \displaystyle {\sum_{\gamma \in \Gamma_2}} \ \lambda (\gamma) \ \displaystyle {\sum_{g \in \gamma}} \ \omega_g (X,Y) \otimes g , \ \forall X, Y, \in V.$ Then $H_{\hbar \lambda} \:= \ (G * T(V)) \ [\hbar] \ / \ {I_\lambda}.$

By definition, the SRA-relations hold in $H_{\hbar \lambda}.$

Theorem (2.16) of [12] ((9.5) of [11]) proves that, when varying
$\lambda$ in $\CC^{\Gamma_2},$ $H_{\hbar \lambda}$ provides an
algebraic deformation of $G * W,$ non trivial, as a deformation, if
$\lambda \neq 0.$ We recall that an algebraic deformation of an
algebra $A$ is a $\CC [\hbar]$-algebra structure on $A [\hbar],$ and
that a polynomial deformation is a formal deformation on $A [[\hbar]]$
such that $A [\hbar]$ is a subalgebra (and therefore an algebraic
deformation), in other words: $\forall a, \ b \in A, \ a *_\hbar b \in
A [\hbar].$ \vskip5pt (4.3) We shall now give a completely different
proof of the E.G. Theorem (2.16) in [12] ((9.5) in [11]), and PBW
Theorem (1.3) in [12] (8.3) in [11]): we show that there exists a $\CC
[G]$-relative deformation of $G * W$, satisfying the SRA-relations,
that, up to an equivalence, this deformation is polynomial, and that
the polynomial part is isomorphic to $H_{\hbar \lambda}$ ; the PBW
property \hskip5pt is a consequence. The proof uses essentially the
C.A. Theorems, classical deformation theory, and a formula of Berezin
([6], [9]).  \vskip5pt (4.4) E.G. THEOREM \vskip5pt {\it There exists
  a non trivial polynomial $\CC [G]$-relative deformation of $G * W$,
  satisfying the SRA-relations. The subalgebra $(G * W)$ [$\hbar$] of
  this deformation is isomorphic to the Symplectic Reflection Algebra
  $H_{\hbar \lambda}$, and the PBW property holds for} $H_{\hbar
  \lambda}$ ($\hbar$ {\it formal, or} $\hbar \in \CC).$ \vskip5pt
$\underline {\rm Proof }$: \vskip5pt (4.4.1) Let $C = C_\lambda.$ As
seen in (4.1), there exists a non-trivial first order $\CC
[G]$-relative deformation: $a \displaystyle\mathop{*}_{\hbar}
\oalign{\cr} b = a * b + \hbar \ C(a,b), \ \forall a,b \in W$ ($*$ is
the Moyal product), and we have to look for a second order still $\CC
[G]$-relative deformation $a \displaystyle\mathop{*}_{\hbar}
\oalign{\cr} b = a * b + \hbar \ C(a,b) + \hbar^2 \ D(a,b),$ such that
$D \ \vert_{\Lambda^2} = 0,$ which is the SRA-condition. Using
Gerstenhaber bracket, the associativity condition at order $2$ is $d D
= - {1 \over 2} \ [C,C].$ Since $C \in Z^2_{\CC [G]} (G * W),$ one has
$[C,C] \in Z^3_{\CC [G]} (G * W) = B^3_{\CC [G]} (G * W),$ so $[C,C] =
dB, \ B \in \Cc^2_{\CC [G]} (G * W),$ moreover, since $C(\Lambda^2)
\subset \CC [G],$ one has $[C,C]_{\vert_{\Lambda^3}} = 0.$ With the
notations of (3.3), and $d_\gamma:= d_{\sigma_\gamma}$ one has
$[\widetilde{C,C}]_\gamma = d_\gamma \tilde{B}_\gamma,
\tilde{B}_\gamma \in \Cc^2_\gamma (W), \ \forall \gamma \in \Gamma
\ \ \ .$ Let $b_\gamma:= \tilde{B}_\gamma \ \vert_{\Lambda^2},$ and
$\Delta_\gamma$ the Koszul differential, then $[\widetilde
  {C,C}]_\gamma \ \vert_{\Lambda^3} = 0$ implies that $d_\gamma
\tilde{B}_\gamma \ \vert_{\Lambda^3} = \Delta_\gamma b_\gamma = 0$, so
$b_\gamma$ is a Koszul cocycle, and it is $S_\gamma$-invariant,
therefore, as in (2.2.12), there exists an $S_\gamma$-invariant
Hochschild cocycle $Z_\gamma \in Z_\gamma^2 (W_\gamma)$ such that
$Z_\gamma \ \vert_{\Lambda^2} = b_\gamma.$ By (3.4.3), there exists $Z
\in Z^2_{\CC [G]} (G * W)$ such that $\tilde{Z}_\gamma = Z_\gamma,
\ \forall \ \gamma \in \Gamma$ ; let $D \:= \ {1 \over 2} \ (Z - B),$
then $dD = - {1 \over 2} [C,C],$ and $\tilde{D}_\gamma
\ \vert_{\Lambda^2} = 0,$ $\forall \ \gamma \in \Gamma,$ moreover $D$
is $\CC [G]$-relative, so $D \ \vert_{\Lambda^2} = 0.$ Now, we have
the wanted second order deformation. The same proof can be repeated
for next orders: for instance, at third order, we have to find $E$
such that $d E = - [D,C],$ and $E \ \vert_{\Lambda^2} = 0.$ Since $D$
is $\CC [G]$-relative, and $D \ \vert_{\Lambda^2} = 0,$ one has $[D,C]
\ \vert_{\Lambda^3} = 0,$ so the above arguments do apply.  \vskip5pt
(4.4.2) We have constructed a non trivial $\CC [G]$-relative
deformation satisfying the SRA-relations. We are going to renormalize,
using an equivalence, to obtain a polynomial deformation still
satisfying SRA-relations. We use arguments inspired of [8] and
[14]. We recall that $*$ is the Moyal product, that . is the abelian
product, and that $W$ is linearly generated by elements of type $X^{*
  k} = X^k, \ X \in V, \ k \in \NN$ (see e.g [19]). Let $\tilde{A} = G
* W \ [[\hbar]],$ with product $*_\hbar$ defined in (4.4.1), there
exists a $\CC [[\hbar]]$-linear map $\rho: W [[\hbar]] \rightarrow
\tilde{A}$ such that: $\rho (1) = 1$ and 
\[\rho (X_1... X_k) = {1 \over {k!}} \ \displaystyle {\sum_{\sigma \in \Sigma_k}} \ X_{\sigma (1)}
\ *_\hbar ... \ *_\hbar \ X_{\sigma (k)}, \ \forall \ X_i \in V,\] and
therefore $\rho (X^k) = X^{{*_\hbar} k}, \ \forall \ X \in V.$ Being
$\CC [G]$-relative, and a deformation of the Moyal product, $*_\hbar$
satisfies:
\[Ad_g (X_1 *_\hbar X_2 *_\hbar... *_\hbar X_k) = g (X_1) *_\hbar g (X_2) *_\hbar... *_\hbar
g (X_k), \forall g \in G, X_i \in V\] (see (3.2.8)), so if we define a
$\CC [[\hbar]]$-linear map $\tilde{\rho}: \tilde{A} \rightarrow
\tilde{A}$ by $\tilde{\rho} (a * g) = \rho (a) * g, \ \forall a \in W,
\ g \in G,$ then $\tilde{\rho} (g * a) = g * \tilde{\rho} (a),$
therefore $\tilde{\rho} = Id + {\sum_k} \ \hbar^k . \rho_k,$ with
$\rho_k \ \CC [G]$-relative. It results that we can define a new $\CC
[G]$-relative product on $\tilde{A}$ by: $\tilde{a} *' \tilde{b} =
\tilde{\rho}^{-1} \ [\tilde{\rho} (\tilde{a}) \ *_\hbar \tilde{\rho}
  (\tilde{b})], \ \tilde{a}, \ \tilde{b} \in \tilde{A}.$ By
definition, one has: $X^{*' k} = X^k = X^{* k}, \ \forall X \in V, \ k
\in \NN,$ and therefore 
\[{1 \over {k!}}  \ \displaystyle {\sum_{\sigma
    \in \Sigma_k}} X_{\sigma (1)} *' X_{\sigma (2)} *' ... *'
X_{\sigma (k)} = X_1 ... X_k\]

for all $X_i \in V$ and $k \in \NN.$ Moreover:
\[[X,Y]_{*'} = \rho^{-1} \ (\omega (X,Y) + \hbar \ C (X \wedge Y)) =
    \omega (X,Y) + \hbar \ C (X \wedge Y),\]

so the SRA-relations are still verified by $*'.$

\vskip5pt (4.4.3) Let us now show that $*'$ is a polynomial
deformation. Let $(W_k, \ k \geq 0)$ by the canonical filtration of
$W.$ One has $X_1 *' X_2 = {1 \over 2} (X_1 *' X_2 + X_2 *' X_1) + {1
  \over 2} \ [X_1, X_2]_{*'} = X_1 . X_2 + {1 \over 2} \ [X_1,
  X_2]_{*'},$ and $[X_1, X_2]_{*'} \in (\CC [G]) [\hbar].$

\vskip5pt (4.4.4) By induction, we assume: $X_1 *' X_2 *' ... *' X_j =
X_1 . X_2 ... X_j + R, \ \forall \ X_i \in V,$ $j \leq k,$ with $R \in
(W_{j-2} \otimes \CC [G]) [\hbar].$ We need the following Lemma (a
direct consequence of a formula of F.A. Berezin [6], [9]):

\vskip5pt (4.4.5) LEMMA 

\vskip5pt {\em Let $A$ be an algebra, $<a_1,...,a_k>:= \frac{1}{k!}
  \sum_{\sigma \in \Sigma_k} a_{\sigma (1)} ... a_{\sigma (k)},$ $a_i
  \in A,$ then:}
\begin{eqnarray*}
& a&.<a_1,...,a_k> = <a, a_1,..., a_k> + \\ & & \sum_{1 \leq j \leq k}
  \frac{B_j}{j!}  \sum_{\underset{\tau \in \Sigma_j}{i_1< ... <i_j}}
  <ad a_{i_{\tau(j)}} ... ad a_{i_{\tau(1)}}(a), a_1,...,
  \hat{a}_{i_1},..., \hat{a}_{i_j},..., a_k >
\end{eqnarray*} {\em where the $B_j$
   are the Bernoulli numbers.}

\vskip5pt $\underline {\rm Proof }$: 

\vskip5pt Consider $A,$ with its natural bracket, as a Lie Algebra,
and let ${\mathcal U}$ be its enveloping algebra. By [9], the formula
is valid in ${\mathcal U}$ (i.e: with product of ${\mathcal U}$). But
there exists an algebra morphism $\mu: {\mathcal U} \rightarrow A$
such that $\mu|_A = Id_A,$ and applying $\mu$ on the formula written
in ${\mathcal U},$ one obtains the formula written in $A,$ as
wanted. \hfill Q.E.D.

\vskip5pt

From the induction assumption (4.4.4), taking first $X_i=X,$ next
ones= $Y,$ and using $X^{*'r}=X^r, \ Y^{*'s} = Y^s,$ it results that:
$X^r \ *' \ Y^s = X^r . Y^s + R',$ with $R' \in (W_{r+s-2} \ \otimes
\ \CC [G]) \ [\hbar],$ if $r+s \leq k,$ and then, since $W_j$ is
linearly generated by $\{X^\ell, \ \ell \leq j, \ X \in V \}:$

\vskip5pt (4.4.6) $a \ *' \ b = a.b + R",$ with $R" \in (W_{r+s-2}
\ \otimes \ \CC [G]) \ [\hbar],$ if $r+s \leq k, \ \forall \ a \in
W_r, \ b \in W_s.$ Now, using Lemma (4.4.5):
\begin{eqnarray*}
X*' (X_1 ... X_k) &=& X \ *' \ < X_1,..., X_k >_{*'} \\ &=& X . X_1
... X_k + \\ && {\sum_{j \leq k}} \ {{B_j} \over {j!}}
\sum_{\underset{\tau \in \Sigma_j}{i_1< ... <i_j}} < ad
\ X_{i_{\tau(j)}} ... ad \ X_{i_{\tau(1)}} (X), X_1,...,
\hat{X}_{i_1},..., \hat{X}_{i_j},..., X_k >_*
\end{eqnarray*}

Since $\xi := ad \ X_{i_{\tau(j)}} ... ad \ X_{i_{\tau(1)}}(X) = ad
\ X_{i_{\tau(j)}} ...  ad \ X_{i_{\tau(2)}} (ad \ X_{i_{\tau(1)}}(X))$
and \linebreak $ad \ X_{i_{\tau(i)}} (X) \in \CC [G] \ [\hbar],$ $\xi$
is a $\CC [\hbar]$-linear combination of terms of type $(V_1 \ *'
... *' \ V_{j-1}) * t,$ with $V_i \in V, \ t \in \CC [G],$ therefore
an element of $(W_{j-1} \otimes \ \CC [G]) [\hbar],$ and the term $ad
\ X_{i_{\tau(j)}} ... ad \ X_{i_{\tau(1)}}(X) *' X_1 *' ... *'
\hat{X}_{i_1} *' ... *' \hat{X}_{i_j} *' ... *' X_k \in (W_{k-1}
\otimes \ \CC [G]) \ [\hbar].$ Similar arguments show that all terms
in the development of 
\[<ad \ X_{i_{\tau(j)}} ... ad\ X_{i_{\tau(1)}}(X),X_1,..., \hat{X}_{i_1},...,
\hat{X}_{i_j},...,X_k>_{*'}\] are elements of $(W_{k-1} \otimes \CC
    [G]) [\hbar].$ Therefore $X \ *' (X_1 ... X_k) = X . X_1 ... X_k +
    S, \ S \in (W_{k-1} \otimes \ \CC [G]) [\hbar].$ Now $X *' X_1 *'
    ... *' X_k = X *' (X_1 ... X_k) + X \ *' R = X . X_1 ... X_k + S +
    X *' R,$ and since $S$ and $X *' R \in (W_{k-1} \otimes \ \CC [G])
    [\hbar],$ our induction is complete.

To conclude the proof, from $X_1 \ *' ... *' \ X_k = X_1 ... X_k + R,$
with $R \in (W_{k-2} \otimes \ \CC [G]) [\hbar], \ \forall \ X_i \in
V, \ k \in \NN,$ we deduce, as was done for (4.4.6), $a \ *' b = ab +
T, \ \forall \ a \in W_r, \ b \in W_s,$ with $T \in (W_{r+s-2} \otimes
\CC [G]) [\hbar].$ This proves that $*'$ is a polynomial deformation.

\vskip5pt
(4.4.7) We now prove that the subalgebra $(W \otimes \CC [G]) [\hbar]$ of $\tilde{A}$ is isomorphic to the Symplectic Reflexion Algebra $H_{\hbar \lambda}$ of (4.2), and the P.B.W. will follow.

We need some notations: $\times_\hbar$ will be the product of
$H_{\hbar \lambda},$ which is generated, as an algebra, by
$\overline{X}, \ \overline{g}$ and $\hbar, \ X \in V, \ g \in G.$ We
denote by $*_\hbar$ the product on $G * W [\hbar]$ constructed (and
denoted by $*'$) in (4.4.2). We denote by $\times$ the product on $G *
W$ coming from the abelian product of $W.$ We also use the notation
$<a_1,...,a_k>$ of (4.4.5). We observe that $R_{\hbar\lambda} (g(X),
g(Y)) = Ad \ g \ (R_{\hbar\lambda}(X,Y)), \ \forall \ X,Y \in V, g \in
G,$ so the natural action of $G$ on $V$ is preserved in the quotient
$H_{\hbar \lambda} = G * T(V) [\hbar] \ / \ {I_\lambda}:
g(\overline{X}) = \overline{g(X)} = \overline{Ad \ g (X)} = Ad
\ \overline{g} (\overline{X}), \ \forall \ g \in G, \ X \in V.$

\vskip5pt (4.4.8) There exists a morphism $\pi: H_{\hbar \lambda}
\mapsto G * W [\hbar] \rightarrow G * W$ (with Moyal product), which
is onto. We define a section $\sigma$ of $\pi$ by:
\[\sigma (X_1 ... X_k \otimes g) = < \overline{X}_1,..., \overline{X}_k>_{\times_\hbar} \times_\hbar \ \overline {g}. \] 
The map $\sigma$ is one to one, so we can identify the spaces $G * W$
and $\sigma (G * W).$ This being done, we have now $G * W \subset
H_{\hbar \lambda}$ and the product $\times$ of $G * W$ becomes:
\vskip5pt (4.4.9) $X_1 ... X_k = < X_1,...,X_k>_{\times_\hbar}
\ \forall \ X_i \in V, \ g \times (X_1...X_k) = (g(X_1) ... g (X_k)
\ \times_\hbar \ g = g \ \times_\hbar \ (X_1 ... X_k),\forall X_i \in
V, \ g \in G.$

Denote by $G * W \ ( \hbar )$ the subspace of $H_{\hbar \lambda}$ of
elements which are polynomial of $\hbar,$ with coefficients in $G *
W,$ by $G * W_k \ ( \hbar )$ elements which are polynomial of $\hbar$
with coefficients in $G * W_k.$ Repeating identically the arguments of
(4.6), we obtain: \vskip5pt (4.4.10) $X_1 \ \times_\hbar \ X_2
\ \times_\hbar... \ \times_\hbar \ X_k = X_1 ... X_k + R,$ with $R \in
G * W_{k-2} (\hbar), \ \forall \ X_i \in V.$ It results that $H_{\hbar
  \lambda} = G * W (\hbar).$ Now, fix any basis $\{ e_1,...,e_{2n} \}$
of $V,$ the natural morphism from $H_{\hbar \lambda}$ onto $G * W
\ [\hbar]$ maps the generator system $\{ e^{i_1}_1 ... e^{in}_n
\times_\hbar g, \ i_1 ... i_n \in \NN, \ g \in G \}$ of the $\CC
      [\hbar]$-algebra $H_{\hbar \lambda}$ onto a basis of $G * W
      [\hbar],$ so we have the wanted isomorphism.  \vskip5pt (4.4.11)
      From the isomorphism $H_{\hbar \lambda} \simeq G * W [\hbar],$
      we deduce that $H_{\hbar \lambda}$ is a deformation of $G * W$,
      and the P.B.W. property for $H_{\hbar \lambda}.$ Given $c \in
      \CC,$ $H_{c \lambda}$ is defined as the quotient of $G * T(V)$
      by relations $R_{c \lambda} (X,Y), \ X,Y \in V$ (see (4.2)). It
      is easy to check that algebras $H_{c \lambda}$ and $H_{\hbar
        \lambda} / H_{\hbar \lambda} (\hbar - c)$ are isomorphic, and
      the P.B.W. property for $H_{c \lambda}$ follows. \hskip160pt
      Q.E.D.  \vfill

\vskip10pt
{\bf References}
\vskip10pt
[ 1 ] Alev J., Farinati M.A, Lambre T., Solotar A.L.: Homologie des invariants d'une alg\`ebre de Weyl sous l'action d'un groupe fini, J. of Algebra, 232 (2000), 564-577.

[ 2 ] Alev J., Lambre T.: Homologie des invariants d'une alg\`ebre de Weyl, $K$-Theory 18 (1999), 401-411.

[ 3 ] Alvarez M.S.: Algebra structure on the Hochschild cohomology of the ring of invariants of a Weyl algebra under a finite group, J. of Algebra 248 (2002), 291-306.

[ 4 ] Arnal D., Ben Amor H., Pinczon G.: The structure of $s \ell (2,1)$-supersymmetry, Pac. J. Math. 165 (1994), 17-49.

[ 5 ] Bayen F., Flato M., Fronsdal C., Lichnerowicz A. Sternheimer D.: Deformation theory and quantization, Ann. Phys. I, II, (1978), 61-110, 111-151.

[ 6 ] Berezin F.A.: Quelques remarques sur l'enveloppe associative d'une alg\`ebre de Lie, Funct. Anal. i evo prilojenie, 1 (1967), 1-14.

[ 7 ] Cartan H., Eilenberg S.: Homological Algebra, Princeton Univ. Press, Princeton NJ, 1956.

[ 8 ] Dito G.: Kontsevich star product on the dual of a Lie algebra, Lett. Math. Phys., 48 (1999), 307-322.

[ 9 ] Dixmier J.: Alg\`ebres Enveloppantes, Gauthier-Villars Paris, 1974.

[ 10 ] Du Cloux F.: Extensions entre repr\'esentations unitaires irr\'eductibles des groupes de Lie nilpotents, Ast\'erisque, 125 (1985), 129-211.

[ 11 ] Etingof P.: Lectures on Calogero-Moser systems, math. QA / 0606233.

[ 12 ] Etingof P., Ginzburg V.: Symplectic reflection algebras, Calogero-Moser space and deformed Harish.Chandra homomorphism. Invent. Math. 147 (2002), 243-348.

[ 13 ] Gestenhaber M., Schack S.D.: Algebraic cohomology and deformation Theory, in: Deformation theory of Algebras and Structures, NATO-ASI Series C.297, Kluwer Academic Publishers, Dordrecht, 1988.

[ 14 ] Gutt S.: An explicit $*$-product on the cotangent bundle of a Lie group, Lett. Math. Phys. 7 (1983), 249-258.

[ 15 ] Montgomery S.: Fixed Rings of Finite Automorphism Groups of Associative Rings, Lect. Notes in Math., vol. 818, Springer-Verlag, New-York / Berlin 1980.

[ 16 ] Nadaud F.: Generalized deformations, Koszul resolutions, Moyal products, Rev. Math. Phys. 10 (5) (1998), 685-704.

[ 17 ] Pinczon G.: Non commutative Deformation Theory, Lett. Math. Phys. 41 (1997), 101-117.

[ 18 ] Pinczon G.: The enveloping algebra of the Lie superalgebra $osp (1,2),$ J. of Algebra 132 (1) (1990), 219-242.

[ 19 ] Pinczon G., Ushirobira R.: Supertrace and Superquadratic Lie structure on the Weyl Algebra, and Applications to Formal Inverse Weyl Transform, Lett. Math. Phys. 74 (2005), 263-291.

[ 20 ] Sridharan R.: Filtered algebras and representations of Lie algebras. Trans. Amer. Math. Soc., 100 (1961), 530-550.

\end{document}